\documentclass[11pt,a4paper,leqno]{amsart}
\usepackage[english]{babel}

\usepackage[latin1]{inputenc}
\usepackage[T1]{fontenc}
\usepackage{amsfonts}
\usepackage{amsmath}
\usepackage{amssymb}
\usepackage{eurosym}
\usepackage{mathrsfs}
\usepackage{palatino}
\usepackage{color}
\usepackage{esint}
\usepackage{url}
\usepackage{verbatim}
\allowdisplaybreaks[4]
\usepackage{enumerate}

\newcommand{\R}{\mathbb{R}}
\newcommand{\C}{\mathbb{C}}

\newcommand{\calC}{\mathcal{C}}

\newcommand{\Z}{\mathbb{Z}}
\newcommand{\E}{\mathbb{E}}

\newcommand{\calR}{\mathcal{R}}

\newcommand{\bla}{\big \langle}
\newcommand{\bra}{\big \rangle}

\numberwithin{equation}{section}

\newcommand{\ud}[0]{\,\mathrm{d}}

\newcommand{\esssup}[0]{\operatornamewithlimits{ess\,sup}}

\newcommand{\BMO}[0]{\operatorname{BMO}}
\newcommand{\bmo}[0]{\operatorname{bmo}}

\newcommand{\loc}[0]{\operatorname{loc}}

\newcommand{\ch}[0]{\operatorname{ch}}

\newcommand{\calD}[0]{\mathcal{D}}

\newcommand{\wt}[1]{{\widetilde{#1}}}

\swapnumbers
\theoremstyle{plain}
\newtheorem{thm}[equation]{Theorem}
\newtheorem{lem}[equation]{Lemma}

\theoremstyle{definition}

\theoremstyle{remark}
\newtheorem{rem}[equation]{Remark}

\author{Kangwei Li}
\address[K.L.]{BCAM (Basque Center for Applied Mathematics), Alameda de Mazarredo 14, 48009 Bilbao, Spain}
\email{kli@bcamath.org}

\author{Henri Martikainen}
\address[H.M.]{Department of Mathematics and Statistics, University of Helsinki, P.O.B. 68, FI-00014 University of Helsinki, Finland}
\email{henri.martikainen@helsinki.fi}

\author{Emil Vuorinen}
\address[E.V.]{Department of Mathematics and Statistics, University of Helsinki, P.O.B. 68, FI-00014 University of Helsinki, Finland}
\email{emil.vuorinen@helsinki.fi}

\title[Quasi--Banach estimates of bilinear bi-parameter commutators]{Quasi--Banach estimates of commutators of bilinear bi-parameter singular integrals: paraproducts}

\makeatletter
\@namedef{subjclassname@2010}{%
  \textup{2010} Mathematics Subject Classification}
\makeatother


\begin{document}

\begin{abstract}
We complete our boundedness theory of commutators of bilinear bi-parameter singular integrals by establishing the following result.
If $T$ is a bilinear bi-parameter singular integral satisfying suitable $T1$ type assumptions, $\|b\|_{\bmo(\R^{n+m})}  = 1$ and
$1 < p, q \le \infty$ and $1/2 < r < \infty$ satisfy $1/p+1/q = 1/r$, then we have
$$
\|[b, T]_1(f_1, f_2)\|_{L^r(\R^{n+m})} \lesssim \|f_1\|_{L^p(\R^{n+m})} \|f_2\|_{L^q(\R^{n+m})}.
$$
Previously the range $r \le 1$ was proved only in the paraproduct free situation. The main novelty lies
in the treatment of the so called partial paraproducts.
\end{abstract}

\maketitle

\section{Introduction}
In this paper we study commutator estimates of general bilinear bi-parameter singular integrals $T$ -- for the precise formulation of such operators, and for additional
background on multilinear multi-parameter analysis, see our previous paper \cite{LMV1}.
There we showed a new dyadic representation theorem under $T1$ type assumptions in the complicated bilinear bi-parameter framework, and used it to conclude various boundedness
properties, including weighted estimates. A model of a bilinear bi-parameter CZO in $\R^n \times \R^m$ is
$$
(T_1 \otimes T_2)(f_1 \otimes f_2, g_1 \otimes g_2)(x) := T_1(f_1, g_1)(x_1)T_2(f_2, g_2)(x_2),
$$
where $f_1, g_1 \colon \R^n \to \C$, $f_2, g_2 \colon \R^m \to \C$,
$x = (x_1, x_2) \in \R^{n+m}$, $T_1$ is a bilinear CZO in $\R^n$ and $T_2$ is a bilinear CZO in $\R^m$. 
For the definition of usual bilinear singular integrals see e.g. Grafakos--Torres \cite{GT}.

We continue our study of commutator estimates in this context, which we initated in \cite{LMV2}.
To this end, define the commutators
$$
[b,T]_1(f_1,f_2) = bT(f_1, f_2) - T(bf_1, f_2) \, \textup{ and } \, [b,T]_2(f_1, f_2) = bT(f_1,f_2) - T(f_1, bf_2). 
$$
A function $b$ is in little BMO, $\bmo(\R^{n+m})$, if $b(\cdot, x_2)$ and $b(x_1, \cdot)$ are uniformly in BMO.
The following is our main theorem.
\begin{thm}\label{thm:main}
Let $1 < p,q \le \infty$ and $1/2 < r < \infty$ satisfy $1/p + 1/q = 1/r$, and let $b \in \bmo(\R^{n+m})$.
Suppose $T$ is a bilinear bi-parameter Calder\'on--Zygmund operator satisfying the assumptions of the bilinear bi-parameter
representation theorem \cite[Theorem 5.1]{LMV1}.
Then we have
$$
\|[b, T]_1(f_1, f_2)\|_{L^r(\R^{n+m})} \lesssim \|b\|_{\bmo(\R^{n+m})} \|f_1\|_{L^p(\R^{n+m})} \|f_2\|_{L^q(\R^{n+m})},
$$
and similarly for $[b, T]_2$.
\end{thm}
\begin{rem}
In \cite{LMV2} we also considered iterated commutators, like $[b_2, [b_1, T]_1]_2$. While we content with first order commutators here, we do not see
any complications in adapting the arguments of \cite{LMV2} and the current paper to prove quasi--Banach estimates for iterated commutators
as well.
\end{rem}

Previously, in \cite{LMV2}, we established this result only for singular integrals that are free of paraproducts (so that they
have a representation using so called bilinear bi-parameter shifts). However, we proved the Banach range boundedness, where $p,q,r \in (1,\infty)$,
for all $T$ as above.
Therefore, the main purpose of this paper is to generalise the result of \cite{LMV2} by removing the paraproduct free
assumption in the quasi--Banach case. This means that we have to deal with the quasi--Banach boundedness
of the commutators of bilinear bi-parameter partial paraproducts and full paraproducts. These are the dyadic model operators, together with shifts,
that appear in the representation theorem of \cite{LMV1}.

The theory of commutator estimates is extremely vast and important. We only mention here some recent related results.
Ou, Petermichl and Strouse proved in \cite{OPS}  that $[b,T] \colon L^2(\R^{n+m}) \to L^2(\R^{n+m})$, when $T$ is a paraproduct
free (linear) bi-parameter singular integral. This was eventually generalised to concern all bi-parameter singular integrals satisfying $T1$ conditions
by Holmes--Petermichl--Wick \cite{HPW} -- in fact, they prove a more general Bloom type two-weight bound. Recently, in \cite{LMV3} we gave
an efficient proof of \cite{HPW} (inspired in part by \cite{LMV2}), and generalised the result of \cite{HPW} to iterated commutators.
All these bi-parameter results rely on the bi-parameter representation theorem \cite{Ma1} by one of us, and on
increasingly sophisticated ways to bound the appearing model operators. Similarly, we crucially rely on our bilinear bi-parameter representation theorem \cite{LMV1} here.
In the one parameter case sparse domination is available as an important tool. For recent one parameter results see e.g.
Holmes--Lacey--Wick \cite{HLW, HLW2},  Lerner--Ombrosi--Rivera-R\'ios \cite{LOR1, LOR2} and Hyt\"onen \cite{Hy5}.

Compared to the linear case one of the key additional difficulties of the bilinear setting lies in obtaining quasi--Banach estimates. These are in general a challenge to obtain in the bi-parameter setting -- even when no commutators are present. This is because the bi-parameter setting removes many standard tools, such as, sparse domination and the Calder\'on--Zygmund decomposition. In this paper we finalise our treatment of bilinear bi-parameter commutator estimates by proving the boundedness in the full range for
all singular integrals satisfying the assumptions of the representation theorem \cite{LMV1}. The main novelty lies in the handling of the partial paraproducts,
which are a tricky class of model operators.

\subsection*{Acknowledgements}
K. Li is supported by Juan de la Cierva - Formaci\'on 2015 FJCI-2015-24547, by the Basque Government through the BERC
2018-2021 program and by Spanish Ministry of Economy and Competitiveness
MINECO through BCAM Severo Ochoa excellence accreditation SEV-2013-0323
and through project MTM2017-82160-C2-1-P funded by (AEI/FEDER, UE) and
acronym ``HAQMEC''.

H. Martikainen is supported by the Academy of Finland through the grants 294840 and 306901, the three-year research grant 75160010 of the University of Helsinki,
and is a member of the Finnish Centre of Excellence in Analysis and Dynamics Research.

E. Vuorinen is supported by the Academy of Finland through the grant 306901 and by the Finnish Centre of Excellence in Analysis and Dynamics Research.

\section{Definitions and preliminaries}
\subsection{Basic notation}
We denote $A \lesssim B$ if $A \le CB$ for some constant $C$ that can depend on the dimension of the underlying spaces, on integration exponents, and on various other constants appearing in the assumptions. We denote $A \sim B$ if $B \lesssim A \lesssim B$.

We work in the bi-parameter setting in the product space $\R^{n+m}$.
In such a context $x = (x_1, x_2)$ with $x_1 \in \R^n$ and $x_2 \in \R^m$.
We often take integral pairings with respect to one of the two variables only:
If $f \colon \R^{n+m} \to \C$ and $h \colon \R^n \to \C$, then $\langle f, h \rangle_1 \colon \R^{m} \to \C$ is defined by
$$
\langle f, h \rangle_1(x_2) = \int_{\R^n} f(y_1, x_2)h(y_1)\ud y_1.
$$

\subsection{Dyadic notation, Haar functions and martingale differences}
We denote a dyadic grid in $\R^n$ by $\calD^n$ and a dyadic grid in $\R^m$ by $\calD^m$. If $I \in \calD^n$, then $I^{(k)}$ denotes the unique dyadic cube $S \in \calD^n$ so that $I \subset S$ and $\ell(S) = 2^k\ell(I)$. Here $\ell(I)$ stands for side length. Also, $\text{ch}(I)$ denotes the dyadic children of $I$ -- this means that $I' \in \ch(I)$ if $(I')^{(1)} = I$. We sometimes write
$\calD = \calD^n \times \calD^m$.

When $I \in \calD^n$ we denote by $h_I$ a cancellative $L^2$ normalised Haar function. This means the following.
Writing $I = I_1 \times \cdots \times I_n$ we can define the Haar function $h_I^{\eta}$, $\eta = (\eta_1, \ldots, \eta_n) \in \{0,1\}^n$, by setting
\begin{displaymath}
h_I^{\eta} = h_{I_1}^{\eta_1} \otimes \cdots \otimes h_{I_n}^{\eta_n}, 
\end{displaymath}
where $h_{I_i}^0 = |I_i|^{-1/2}1_{I_i}$ and $h_{I_i}^1 = |I_i|^{-1/2}(1_{I_{i, l}} - 1_{I_{i, r}})$ for every $i = 1, \ldots, n$. Here $I_{i,l}$ and $I_{i,r}$ are the left and right
halves of the interval $I_i$ respectively. If $\eta \in \{0,1\}^n \setminus \{0\}$ the Haar function is cancellative: $\int h_I^{\eta} = 0$. We usually suppress the presence of $\eta$
and simply write $h_I$ for some $h_I^{\eta}$, $\eta \in \{0,1\}^n \setminus \{0\}$. Then $h_Ih_I$ can stand for $h_I^{\eta_1} h_I^{\eta_2}$, but we always treat
such a product as a non-cancellative function (which it is in the worst case scenario $\eta_1 = \eta_2$).

For $I \in \calD^n$ and a locally integrable function $f\colon \R^n \to \C$, we define the martingale difference
$$
\Delta_I f = \sum_{I' \in \textup{ch}(I)} \big[ \bla f \bra_{I'} -  \bla f \bra_{I} \big] 1_{I'}.
$$
Here $\bla f \bra_I = \frac{1}{|I|} \int_I f$. We also write $E_I f = \bla f \bra_I 1_I$.
Now, we have $\Delta_I f = \sum_{\eta \ne 0} \langle f, h_{I}^{\eta}\rangle h_{I}^{\eta}$, or suppressing the $\eta$ summation, $\Delta_I f = \langle f, h_I \rangle h_I$, where $\langle f, h_I \rangle = \int f h_I$. A martingale block is defined by
$$
\Delta_{K,i} f = \mathop{\sum_{I \in \calD^n}}_{I^{(i)} = K} \Delta_I f, \qquad K \in \calD^n.
$$

Next, we define bi-parameter martingale differences. Let $f \colon \R^n \times \R^m \to \C$ be locally integrable.
Let $I \in \calD^n$ and $J \in \calD^m$. We define the martingale difference
$$
\Delta_I^1 f \colon \R^{n+m} \to \C, \Delta_I^1 f(x) := \Delta_I (f(\cdot, x_2))(x_1).
$$
Define $\Delta_J^2f$ analogously, and also define $E_I^1$ and $E_J^2$ similarly.
We set
$$
\Delta_{I \times J} f \colon \R^{n+m} \to \C, \Delta_{I \times J} f(x) = \Delta_I^1(\Delta_J^2 f)(x) = \Delta_J^2 ( \Delta_I^1 f)(x).
$$
Notice that $\Delta^1_I f = h_I \otimes \langle f , h_I \rangle_1$, $\Delta^2_J f = \langle f, h_J \rangle_2 \otimes h_J$ and
$ \Delta_{I \times J} f = \langle f, h_I \otimes h_J\rangle h_I \otimes h_J$ (suppressing the finite $\eta$ summations).
Martingale blocks are defined in the natural way
$$
\Delta_{K \times V}^{i, j} f  =  \sum_{I\colon I^{(i)} = K} \sum_{J\colon J^{(j)} = V} \Delta_{I \times J} f = \Delta_{K,i}^1( \Delta_{V,j}^2 f) = \Delta_{V,j}^2 ( \Delta_{K,i}^1 f).
$$

\subsection{Weights}
We need weights for certain things, even though they do not appear in the main theorem.
A weight $w(x_1, x_2)$ (i.e. a locally integrable a.e. positive function) belongs to the bi-parameter $A_p$ class, $A_p(\R^n \times \R^m)$, $1 < p < \infty$, if
$$
[w]_{A_p(\R^n \times \R^m)} := \sup_{R} \bla w \bra_R \bla w^{1-p'} \bra_R^{p-1} < \infty,
$$
where the supremum is taken over $R = I \times J$, where $I \subset \R^n$ and $J \subset \R^m$ are cubes
with sides parallel to the axes (we simply call such $R$ rectangles).
We have
$$
[w]_{A_p(\R^n\times \R^m)} < \infty \textup { iff } \max\big( \esssup_{x_1 \in \R^n} \,[w(x_1, \cdot)]_{A_p(\R^m)}, \esssup_{x_2 \in \R^m}\, [w(\cdot, x_2)]_{A_p(\R^n)} \big) < \infty,
$$
and that $\max\big( \esssup_{x_1 \in \R^n} \,[w(x_1, \cdot)]_{A_p(\R^m)}, \esssup_{x_2 \in \R^m}\, [w(\cdot, x_2)]_{A_p(\R^n)} \big) \le [w]_{A_p(\R^n\times \R^m)}$, while
the constant $[w]_{A_p}$ is dominated by the maximum to some power.
Of course, $A_p(\R^n)$ is defined similarly as $A_p(\R^n \times \R^m)$ -- just take the supremum over cubes $Q$. For the basic theory
of bi-parameter weights consult e.g. \cite{HPW}.

We record the following standard weighted square function estimates.
\begin{lem}\label{lem:standardEst1}
For $p \in (1,\infty)$ and $w \in A_p(\R^n \times \R^m)$ we have
\begin{align*}
\| f \|_{L^p(w)}
& \sim_{[w]_{A_p(\R^n \times \R^m)}} \Big\| \Big( \mathop{\sum_{I \in \calD^n}}_{J \in \calD^m} |\Delta_{I \times J} f|^2 \Big)^{1/2} \Big\|_{L^p(w)} \\
&\sim_{[w]_{A_p(\R^n \times \R^m)}}  \Big\| \Big(  \sum_{I \in \calD^n} |\Delta_I^1 f|^2 \Big)^{1/2} \Big\|_{L^p(w)}
\sim_{[w]_{A_p(\R^n \times \R^m)}} \Big\| \Big(  \sum_{J \in \calD^m} |\Delta_J^2 f|^2 \Big)^{1/2} \Big\|_{L^p(w)}.
\end{align*}
\end{lem}

\subsection{Maximal functions}
Given $f \colon \R^{n+m} \to \C$ and $g \colon \R^n \to \C$ we denote the dyadic maximal functions
by
$$
M_{\calD^n}g(x) := \sup_{I \in \calD^n} \frac{1_I(x)}{|I|}\int_I |g(y)| \ud y
$$
and
$$
M_{\calD^n, \calD^m} f(x_1, x_2) := \sup_{R \in \calD^n \times \calD^m}  \frac{1_R(x_1, x_2)}{|R|}\iint_R |f(y_1, y_2)|\ud y_1 \ud y_2.
$$
The non-dyadic variants are simply denoted by $M$, as it is clear what is meant from the context.
We also set $M^1_{\calD^n} f(x_1, x_2) =  M_{\calD^n}(f(\cdot, x_2))(x_1)$. The operator $M^2_{\calD^m}$ is defined similarly.

Standard weighted estimates involving maximal functions are recorded in the following lemma.
\begin{lem}\label{lem:standardEst2}
For $p, s \in (1,\infty)$ and $w \in A_p$ we have the Fefferman--Stein inequality
$$
\Big\| \Big( \sum_j |M f_j |^s \Big)^{1/s} \Big\|_{L^p(w)} \le C([w]_{A_p}) \Big\| \Big( \sum_{j} | f_j |^s \Big)^{1/s} \Big\|_{L^p(w)}.
$$
We also have
$$
\| \varphi_{\calD^n}^1 f\|_{L^p(w)} 
\sim_{[w]_{A_p}}
\Big\| \Big( \sum_{I \in \calD^n} \frac{1_I}{|I|} \otimes [M \langle f, h_I \rangle_1]^2 \Big)^{1/2} \Big\|_{L^p(w)} 
\lesssim_{[w]_{A_p}} \|f\|_{L^p(w)},
$$
where
$$
\varphi_{\calD^n}^1 f  := \sum_{I \in \calD^n} h_I \otimes M \langle f, h_I \rangle_1.
$$
The function $\varphi_{\calD^m}^2 f$ is defined in the symmetric way and satisfies the same estimates.
\end{lem}

\subsection{BMO spaces}\label{ss:bmo}
We say that a locally integrable function $b \colon \R^n \to \C$ belongs to the dyadic BMO space $\BMO(\calD^n)$ if
$$
\|b\|_{\BMO(\calD^n)} := \sup_{I \in \calD^n} \frac{1}{|I|} \int_I |b - \langle b \rangle_I| < \infty.
$$
The ordinary space $\BMO(\R^n)$ is defined by taking the supremum over all cubes. 

We say that a locally integrable function $b \colon \R^{n+m} \to \C$ belongs to the dyadic little BMO space $\bmo(\calD)$,
where $\calD = \calD^n \times \calD^m$, if
$$
\|b\|_{\bmo(\calD)} := \sup_{R \in \calD} \frac{1}{|R|} \int_R |b - \langle b \rangle_R| < \infty.
$$
The non-dyadic space $\bmo(\R^{n+m})$ is defined in the natural way -- take the supremum over all rectangles. We have
$$
\|b\|_{\bmo(\calD)} \sim \max\big( \esssup_{x_1 \in \R^n} \, \|b(x_1, \cdot)\|_{\BMO(\calD^m)}, \esssup_{x_2 \in \R^m}\, \|b(\cdot, x_2)\|_{\BMO(\calD^n)} \big).
$$

Finally, we have the product BMO space. Set
$$
\|b\|_{\BMO_{\textup{prod}}(\calD)} := 
\sup_{\Omega} \Big( \frac{1}{|\Omega|} \mathop{\sum_{I \times J \in \calD}}_{I \times J \subset \Omega} |\langle b, h_I \otimes h_J\rangle|^2  \Big)^{1/2},
$$
where the supremum is taken over those sets $\Omega \subset \R^{n+m}$ such that $|\Omega| < \infty$ and such that for every $x \in \Omega$ there exist
$I \times J \in \calD$ so that $x \in I \times J \subset \Omega$.
The non-dyadic product BMO space $\BMO_{\textup{prod}}(\R^{n+m})$ can be defined using the norm defined by the supremum over all dyadic grids of
the above dyadic norms. 

\subsection{Adapted maximal functions}\label{sec:AdapMaxFunc}
For $b \in \BMO(\R^n)$ and $f \colon \R^n \to \C$ define
$$
M_bf = \sup_I \frac{1_I}{|I|} \int_I |b-\langle b \rangle_I| |f|.
$$
In the situation $b \in \bmo(\R^{n+m})$ and $f \colon \R^{n+m} \to \C$ we similarly define
$$
M_b f = \sup_{I,J} \frac{1_{I \times J}}{|I||J|} \iint_{I \times J} |b-\langle b \rangle_{I \times J}| |f|.
$$
Here the supremums are taken over all intervals $I \subset \R^n$ and $J \subset \R^m$. The dyadic variants could also be defined, and denoted by
$M_{\calD^n, b}$ and $M_{\calD^n, \calD^m, b}$.

For a little BMO function $b \in \bmo(\R^{n+m})$ define 
$$
\varphi_{\calD^m, b}^2(f) = \sum_{J \in \calD^m} M_{\langle b \rangle_{J,2}} \langle f, h_J \rangle_2 \otimes h_J,
$$
and similarly define $\varphi_{\calD^n, b}^1(f)$. For our later usage it is important to not to use the dyadic
variant $M_{\calD^n, \langle b \rangle_{J,2}}$, as it would induce an unwanted dependence on $\calD^n$ (which has relevance
in some randomisation considerations). For the following lemma see \cite{LMV2}.
\begin{lem}\label{lem:bmaxbounds}
Suppose $\|b_i\|_{\BMO(\R^n)} \le 1$, $1 < s, p < \infty$ and $w \in A_p(\R^n)$. Then we have
\begin{equation}\label{eq:vMb}
\Big\| \Big( \sum_i [M_{b_i} f_i]^s \Big)^{1/s} \Big\|_{L^p(w)} \lesssim C([w]_{A_p(\R^n)}) \Big\| \Big( \sum_i |f_i|^s \Big)^{1/s} \Big\|_{L^p(w)}.
\end{equation}
The same bound holds with $\|b_i\|_{\bmo(\R^n \times \R^m)} \le 1$ and $w \in A_p(\R^n \times \R^m)$.
For a function $b$ with $\|b\|_{\bmo(\R^n \times \R^m)} \le 1$ we also have
$$
\|\varphi_{\calD^m, b}^2(f)\|_{L^p(w)} \le C([w]_{A_p(\R^n \times \R^m)})\|f\|_{L^p(w)}, \qquad 1 < p < \infty,\, w \in A_p(\R^n \times \R^m).
$$
\end{lem}

\subsection{Commutators}
We set 
$$
[b,T]_1(f_1,f_2) = bT(f_1, f_2) - T(bf_1, f_2) \, \textup{ and } \, [b,T]_2(f_1, f_2) = bT(f_1,f_2) - T(f_1, bf_2). 
$$
These are understood in a situation, where we e.g. already know that $T \colon L^3(\R^{n+m}) \times L^3(\R^{n+m}) \to L^{3/2}(\R^{n+m})$, and
$b$ is locally in $L^3$. Then we initially study the case that $f_1$ and $f_2$ are, say, bounded and compactly supported, so that
e.g. $bf_2 \in L^3(\R^{n+m})$ and $bT(f_1,f_2) \in L^1_{\loc}(\R^{n+m})$.

\subsection{Random dyadic grids}
We need the following notation regarding random dyadic grids.
Let $\mathcal{D}_0^n$ and $\mathcal{D}_0^m$ denote the standard dyadic grids on $\R^n$ and $\R^m$ respectively.
For $\omega_1 = (\omega^i_1)_{i \in \Z} \in (\{0,1\}^n)^{\Z}$, $\omega_2 = (\omega^i_2)_{i \in \Z} \in(\{0,1\}^m)^{\Z}$, $I \in \calD^n_0$ and $J \in \calD^m_0$ denote
$$
I + \omega_1 := I + \sum_{i:\, 2^{-i} < \ell(I)} 2^{-i}\omega_1^i \qquad \textup{and} \qquad J + \omega_2 := J + \sum_{i:\, 2^{-i} < \ell(J)} 2^{-i}\omega_2^i.
$$
Then we define the random lattices
$$\calD^n_{\omega_1} = \{I + \omega_1\colon I \in \calD^n_0\} \qquad \textup{and} \qquad
\calD^m_{\omega_2} = \{J + \omega_2\colon J \in \calD^m_0\}.
$$
There is a natural probability product measure $\mathbb{P}_{\omega_1}$ in $(\{0,1\}^n)^{\Z}$ and $\mathbb{P}_{\omega_2}$ in $(\{0,1\}^m)^{\Z}$.
We set $\omega = (\omega_1, \omega_2) \in (\{0,1\}^n)^{\Z} \times (\{0,1\}^m)^{\Z}$,
and denote the expectation over the product probability space by $\E_{\omega} = \E_{\omega_1, \omega_2} = \E_{\omega_1} \E_{\omega_2} = \E_{\omega_2} \E_{\omega_1} =
\iint \ud \mathbb{P}_{\omega_1} \ud \mathbb{P}_{\omega_2}$.

We also set $\calD_0 = \calD^n_0 \times \calD^m_0$. Given $\omega = (\omega_1, \omega_2)$ and $R = I \times J \in \calD_0$ we may set
$$
R + \omega = (I+\omega_1) \times (J+\omega_2) \qquad \textup{and} \qquad \calD_{\omega} = \{R + \omega\colon \, R \in \calD_0\}.
$$
\section{Martingale difference expansions of products and a duality lemma}\label{sec:marprod}
\subsection{Expansions of products}
We begin by recalling from \cite{LMV2} our strategy of expanding biparameter commutators of dyadic model operators, which differs from the strategy used in \cite{HPW} in some key ways.
A product $bf$ paired with Haar functions is expanded in the bi-parameter fashion only if both of the Haar functions are cancellative. In a mixed
situation we expand only in $\R^n$ or $\R^m$, and in the remaining fully non-cancellative situation we do not expand at all.
Our protocol also entails the following: when pairing with a non-cancellative Haar function we add and subtract a suitable average of $b$.

Let $\calD^n$ and $\calD^m$ be some fixed dyadic grids in $\R^n$ and $\R^m$, respectively, and write $\calD= \calD^n \times \calD^m$.
In what follows we sum over $I \in \calD^n$ and $J \in \calD^m$.

\subsubsection*{Paraproduct operators}
We define certain standard paraproduct operators:
\begin{align*}
A_1(b,f) &= \sum_{I, J} \Delta_{I \times J} b \Delta_{I \times J} f, \,\,
A_2(b,f) = \sum_{I, J} \Delta_{I \times J} b E_I^1\Delta_J^2 f, \\
A_3(b,f) &= \sum_{I, J} \Delta_{I \times J} b \Delta_I^1 E_J^2  f, \,\,
A_4(b,f) = \sum_{I, J} \Delta_{I \times J} b \bla f \bra_{I \times J},
\end{align*}
and
\begin{align*}
A_5(b,f) &= \sum_{I, J} E_I^1 \Delta_J^2 b \Delta_{I \times J} f, \,\,
A_6(b,f) = \sum_{I, J}  E_I^1 \Delta_J^2 b  \Delta_I^1 E_J^2  f, \\
A_7(b,f) &= \sum_{I, J} \Delta_I^1 E_J^2  b \Delta_{I \times J} f, \,\,
A_8(b,f) = \sum_{I, J}  \Delta_I^1 E_J^2 b E_I^1 \Delta_J^2 f.
\end{align*}
The operators are grouped into two collections, since they are handled differently (using product BMO or little BMO estimates, respectively).
Also recall that $\bmo \subset \BMO_{\textup{prod}}$.
When the underlying grid needs to be written, we write $A_{i, \calD}(b, f)$.

We also define
$$
a^1_1(b,f) = \sum_I \Delta_I^1 b \Delta_I^1 f \qquad \textup{and} \qquad
a^1_2(b,f) = \sum_I \Delta_I^1 b E_I^1 f.
$$
The operators $a^2_1(b,f)$ and $a^2_2(b,f)$ are defined analogously. Again, we can also e.g. write
$a^1_{1, \calD^n}(b,f)$ to emphasise the underlying dyadic grid.

For the following standard lemma see \cite{HPW} or \cite{LMV2} (the operators
$a_j^1(b,\cdot)$, $a_j^2(b,\cdot)$ do not appear in \cite{HPW}).
\begin{lem}\label{lem:basicAa}
Let $H_b$ be $A_i(b,\cdot)$, $i=1,\ldots, 8$, or $a_j^1(b,\cdot)$, $a_j^2(b,\cdot)$, $j=1,2$. Let also $b\in \bmo(\R^{n+m})$, $p \in (1,\infty)$ and $w \in A_p(\R^n \times \R^m)$.
Then we have
\begin{align*}
\|H_b f\|_{L^p(w)} \lesssim C([w]_{A_p(\R^n \times \R^m)}) \|b\|_{\bmo(\R^{n+m})} \|f\|_{L^p(w)}.
\end{align*}
\end{lem}
Let now $f \in L^p(\R^{n+m})$ for some $p \in (1,\infty)$, and $b \in \bmo(\R^{n+m})$.
We know that $b \in L^{p}_{\loc}(\R^{n+m})$ by the John--Nirenberg valid for little BMO (see \cite{HPW}).
For $I_0 \in \calD^n$ and $ J_0 \in \calD^m$ we will now introduce our expansions of $\langle bf, h_{I_0} \otimes h_{J_0}\rangle$, $\big \langle bf, h_{I_0} \otimes \frac{1_{J_0}}{|J_0|}\big\rangle$ and $\langle bf \rangle_{I_0 \times J_0}$.
\subsubsection*{Expansion of $\langle bf, h_{I_0} \times h_{J_0} \rangle$}
There holds
$$
1_{I_0 \times J_0} b
= \sum_{\substack{I_1\times J_1 \in \calD \\ I_1 \times J_1 \subset I_0 \times J_0}}\Delta_{I_1 \times J_1} b
+\sum_{\substack{J_1 \in \calD^m \\ J_1 \subset J_0}} E^1_{I_0} \Delta^2_{J_1} b
+ \sum_{\substack{I_1 \in \calD^n \\ I_1 \subset I_0}} \Delta^1_{I_1} E^2_{J_0} b
+ E_{I_0 \times J_0} b.
$$
Let us denote these terms by $I_j$, $j=1,2,3,4$, in the respective order.
We have the corresponding decomposition of $f$, whose terms we denote by $II_i$, $i=1,2,3,4$. Calculating carefully the pairings
$\langle I_j II_i, h_{I_0} \otimes h_{J_0} \rangle$ we see that
\begin{equation}\label{eq:biparEX}
\langle bf, h_{I_0} \otimes h_{J_0} \rangle = \sum_{i=1}^8 \langle A_i(b, f), h_{I_0} \otimes h_{J_0} \rangle + \langle b \rangle_{I_0 \times J_0} \langle f, h_{I_0} \otimes h_{J_0} \rangle.
\end{equation}
\subsubsection*{Expansion of $\big \langle bf, h_{I_0} \otimes \frac{1_{J_0}}{|J_0|}\big\rangle$}
This time we write
$
1_{I_0} b  = \sum_{\substack{I_1 \in \calD^n \\ I_1 \subset I_0}}\Delta_{I_1}^1 b + E_{I_0}^1 b,
$
and similarly for $f$. Calculating $\langle bf, h_{I_0} \rangle_1$ we see that
\begin{equation}\label{eq:1EX}
\begin{split}
\Big \langle bf, h_{I_0} \otimes \frac{1_{J_0}}{|J_0|}\Big\rangle &= \sum_{i=1}^2 \Big\langle a_i^1(b,f), h_{I_0} \otimes \frac{1_{J_0}}{|J_0|} \Big\rangle \\
&+ \bla (\langle b \rangle_{I_0,1} - \langle b \rangle_{I_0 \times J_0}) \langle f, h_{I_0}\rangle_1\bra_{J_0}
+ \langle b \rangle_{I_0 \times J_0} \Big \langle f, h_{I_0} \otimes \frac{1_{J_0}}{|J_0|}\Big\rangle.
\end{split}
\end{equation}

When we have $\langle bf \rangle_{I_0 \times J_0}$ we do not expand at all:
\begin{equation}\label{eq:noEX}
\langle bf \rangle_{I_0 \times J_0} = \langle (b-\langle b \rangle_{I_0 \times J_0})f \rangle_{I_0 \times J_0}
+ \langle b \rangle_{I_0 \times J_0} \langle f\rangle_{I_0 \times J_0}.
\end{equation}
All of our commutators are simply decomposed using \eqref{eq:biparEX}, \eqref{eq:1EX} (and its symmetric form) and \eqref{eq:noEX} whenever
the relevant pairings/averages appear.

\subsection{Duality lemma}
We present a technical modification of the standard $H^1$-$\BMO$ type duality estimate. A special case of this, formulated below in Remark
\ref{lem:H1-BMO-Modified-1par}, is key when we deal with partial paraproducts.

For a sequence of scalars $\{a_{R}\}_{R \in \calD_{\omega}}$ denote
$$
\|\{a_{R}\}_{R \in \calD_{\omega}}\|_{\BMO_{\textup{prod}}(\calD_{\omega})} := 
\sup_{\Omega} \Big( \frac{1}{|\Omega|} \mathop{\sum_{R \in \calD_{\omega}}}_{R \subset \Omega} |a_{R}|^2  \Big)^{1/2},
$$
where the supremum is taken over those sets $\Omega \subset \R^{n+m}$ such that $|\Omega| < \infty$ and such that for every $x \in \Omega$ there exist
$K \times V \in \calD_{\omega}$ so that $x \in K \times V \subset \Omega$.
\begin{lem}\label{lem:H1-BMO-Modified}
Let $\omega=(\omega_1,\omega_2)$ be a random parameter and $F \subset \R^{n+m}$. 
Suppose $\calC \subset \calD_0$ is a collection of rectangles such that
$|R \cap F | \ge \frac{99}{100} |R|$ for all $R \in \calC$.
Let  $\{a_{R + \omega }\}_{R \in \calC}$ and $\{b_{R}\}_{R \in \calC}$ be two collections of scalars.
Then
$$
\sum_{R \in \calC} | a_{R + \omega} b_{R} |
\lesssim \| \{ a_{R+ \omega} \}_{R \in \calC} \|_{\BMO_{\textup{prod}}(\calD_\omega)}
\iint_F \Big( \sum_{R \in \calC} |b_{R}|^2 \frac{1_{R}}{|R|} \Big)^{1/2}.
$$
\end{lem}

\begin{proof}
Write
$$
S = \Big( \sum_{R \in \calC}|b_{R}|^2 \frac{1_{R}}{|R|} \Big)^{1/2}.
$$
We may suppose that $\| 1_F S \|_{L^1} < \infty$.
For $u \in \Z$ let $\Omega_u=\{ 1_F S> 2^{-u}\}$ and $\wt{\Omega}_u= \{M1_{\Omega_u} >c\}$,
where $c=c(n,m) \in (0,1)$ is small enough. 
Define the collections
$$
\widehat{\calC}_u
=\Big\{ R \in \calC \colon |R \cap \Omega_u| \ge \frac{1}{100} |R|\Big\}
$$
and write $\calC_u=\widehat{\calC}_u \setminus \widehat{\calC}_{u-1}$. 

Let $R \in \calC$ be such that $b_{R} \not=0$. Then for all
$x=(x_1,x_2) \in R \cap F$ there holds that
$$
0 < \frac{|b_{R}| }{|R|^{1/2}} \le 1_F(x) S(x).
$$
Since $|R \cap F | \ge \frac{99}{100} |R|$, this implies that  $R \in \widehat{\calC}_u$ for all large  $u \in \Z$.
On the other hand, since $|\Omega_u| \to 0$, as $u \to -\infty$, we have $R \not \in \widehat{\calC}_u$ for all small $u \in \Z$.
Therefore, (because $\Omega_{u-1} \subset \Omega_u$ for all $u$) it holds that
$$
\sum_{R \in \calC} | a_{R+\omega} b_{R} |
= \sum_{u \in \Z} \sum_{R \in \calC_u} | a_{R+\omega} b_{R} |.
$$

For fixed $u$ estimate
$$
\sum_{R \in \calC_u} | a_{R+\omega} b_{R} |
\le \Big( \sum_{R \in \calC_u}  |a_{R+\omega}|^2\Big)^{1/2}
\Big(\sum_{R \in \calC_u} |b_{R}| ^2\Big)^{1/2}. 
$$
Suppose $R \in \calC_u$. Because $|R \cap \Omega_u| \ge \frac{1}{100}|R|$, there holds that 
$
R+\omega \subset 3R \subset \wt \Omega_u
$
if $c = c(n,m)$ is fixed small enough.
Therefore, since $|\wt \Omega_u| \lesssim |\Omega_u|$, we get that
$$
\Big( \sum_{R \in \calC_u}  |a_{R+\omega}|^2\Big)^{1/2}
\lesssim  \| \{ a_{R+ \omega} \}_{R \in \calC} \|_{\BMO_{\textup{prod}}(\calD_\omega)} | \Omega_u |^{1/2}.
$$

Let again $R \in \calC_u$. Every $R \in \calC$ satisfies by assumption that $|R \cap F| \ge \frac{99}{100} |R|$, and because $R \not \in \widehat \calC_{u-1}$,
there holds that $|R \cap \Omega_{u-1}^c| \ge \frac{99}{100} |R|$. Thus, we have $|R \cap F \cap \Omega_{u-1}^c| \ge \frac{98}{100} |R|$, and this gives
(noting also that $R \subset \wt \Omega_u$ for every $R \in \calC_u$) that
$$
\Big(\sum_{R \in \calC_u} |b_{R}|^2\Big)^{1/2}
\lesssim \Big(\iint \displaylimits _{\wt \Omega_u \setminus \Omega_{u-1}}1_F\sum_{R \in \calC_u} |b_{R}| ^2\frac{1_R}{|R|}\Big)^{1/2}
\le \Big(\iint \displaylimits _{\wt \Omega_u \setminus \Omega_{u-1}}1_FS^2\Big)^{1/2} 
\lesssim 2^{-u}|\Omega_u|^{1/2}.
$$
The claim follows by summing over $u$.
\end{proof}

\begin{rem}\label{lem:H1-BMO-Modified-1par}
Set
$$
\|\{a_{V}\}_{V \in \calD^m_{\omega_2}} \|_{\BMO(\calD^m_{\omega_2})} := 
\sup_{V_0 \in \calD^m_{\omega_2}}\Big( \frac{1}{|V_0|} \sum_{\substack{V \in \calD^m_{\omega_2} \\ V \subset V_0}} |a_V|^2\Big)^{1/2}.
$$

Related to weak type estimates of partial paraproducts we use the following special case of Lemma \ref{lem:H1-BMO-Modified}.
Let $F \subset \R^{n+m}$ and $K_0 \in \calD^n_0$. Suppose $\calC \subset \calD^m_0$ is a collection of cubes such that
$|(K_0 \times V) \cap F | \ge \frac{99}{100} |K_0 \times V|$ for every $V \in \calC$. Let $\omega_2$ be a random parameter,
and let $\{a_{V+\omega_2}\}_{V \in \calC}$ be a collection of scalars. Then, for all scalars  $\{b_V\}_{V \in \calC}$
we have
$$
\sum_{V \in \calC} | a_{V+\omega_2}b_V| 
\lesssim \| \{ a_{V+\omega_2} \}_{V \in \calC} \|_{\BMO(\calD^m_{\omega_2})} \iint_F \frac{1_{K_0}}{|K_0|} \otimes \Big( \sum_{V \in \calC}  |b_V|^2 \frac{1_V}{|V|} \Big)^{1/2}.
$$
This follows from Lemma \ref{lem:H1-BMO-Modified}. To see this, let $\wt \calC=\{K \times V \in \calD_0  \colon K=K_0,  V \in \calC\}$.
For $K \times V \in \calC$ define $\wt a_{K \times (V+\omega_2)} = a_{V+\omega_2}$. 
Then
$$
\| \{\wt a_{K \times (V+\omega_2)}\}_{K \times V \in \wt \calC} \|_{\BMO_{\textup{prod}}(\calD^n_0 \times \calD^m_{\omega_2})}
= \frac{1}{|K_0|^{1/2}} \| \{a_{V+\omega_2}\}_{V \in \calC} \|_{\BMO( \calD^m_{\omega_2})},
$$
from which the claim follows.
\end{rem}

\section{Commutators of partial paraproducts}\label{sec:ComPartialP}
We begin by defining partial paraproducts on random dyadic grids. Let $\omega = (\omega_1, \omega_2) \in (\{0,1\}^n)^{\Z} \times (\{0,1\}^m)^{\Z}$ and
$k = (k_1, k_2, k_3)$, $k_1, k_2, k_3 \ge 0$. For each $K, I_1, I_2, I_3 \in \calD^n_{\omega_1}$ with $I_i^{(k_i)} = K$ and $V \in \calD^m_{\omega_2}$
we are given a constant $a_{K, V, (I_i)}^{\omega}$ such that for all $K, I_1, I_2, I_3 \in \calD^n_{\omega_1}$ with $I_i^{(k_i)} = K$ we have
$$
\|\{ a_{K, V, (I_i)}^{\omega} \}_{V \in \calD^m_{\omega_2}} \|_{\BMO(\calD^m_{\omega_2})}
\le \frac{|I_1|^{1/2} |I_2|^{1/2} |I_3|^{1/2}}{|K|^2}.
$$
A partial paraproduct $P_{k, \omega}$ of complexity $k$ of a particular form is
\begin{equation}\label{eq:ParticularPartialP}
\begin{split}
\langle P_{k, \omega}(f_1,f_2), f_3\rangle =  \sum_{K \in \calD^n_{\omega_1}} \sum_{\substack{I_1, I_2, I_3 \in \calD^n_{\omega_1} \\ I_i^{(k_i)} = K}} \sum_{V \in \calD^m_{\omega_2}}
&a_{K, V, (I_i)}^{\omega} \bla  f_1, h_{I_1}^0 \otimes h_V  \bra  \\
&\times \Big\langle f_2, h_{I_2} \otimes \frac{1_V}{|V|} \Big\rangle \Big\langle f_3, h_{I_3} \otimes \frac{1_V}{|V|} \Big\rangle.
\end{split}
\end{equation}An operator of the above form, but formed using $(h_{I_1}, h_{I_2}^0, h_{I_3})$ or $(h_{I_1}, h_{I_2}, h_{I_3}^0)$
instead of $(h_{I_1}^0, h_{I_2}, h_{I_3})$, or formed using $(1_V/|V|, h_V, 1_V/|V|)$ or $(1_V/|V|, 1_V/|V|, h_V)$ instead of
$(h_V, 1_V/|V|, 1_V/|V|)$, is also a partial paraproduct. So there are nine different types of partial paraproducts. 
Of course, we also have the symmetric partial paraproducts with the shift structure in $\R^m$ and the paraproduct structure in $\R^n$.

We are ready to state our result concerning commutators of partial paraproducts. For technical reasons that appear later
when we want to move our bounds from model operators to singular integrals using the representation \cite{LMV1}
it is necessary to consider averages of random partial paraproducts.
\begin{thm}\label{thm:com1PP}
Let $\|b\|_{\bmo(\R^{n+m})}  = 1$, and let $1 < p, q \le \infty$ and $1/2 < r < \infty$ satisfy $1/p+1/q = 1/r$.
Suppose that $(P_{k, \omega})_{\omega}$ is a collection of partial paraproducts of the same type and of fixed complexity $k = (k_1, k_2, k_3)$.
Then we have
$$
 \|\mathbb{E}_{\omega}[b,P_{k,\omega}]_1(f_1, f_2)\|_{L^r(\R^{n+m})} 
 \lesssim (1+\max k_i) \|f_1\|_{L^p(\R^{n+m})} \|f_2\|_{L^q(\R^{n+m})}.
$$
\end{thm}

\begin{proof}
The case $p,q,r \in (1,\infty)$ in Theorem \ref{thm:com1PP} is easy, since we already know 
the Banach range boundedness of commutators of partial paraproducts by \cite{LMV2}.
Indeed, if $f_1 \in L^p(\R^{n+m})$, $f_2 \in L^q(\R^{n+m})$ and $f_3 \in L^{r'}(\R^{n+m})$, then by \cite{LMV2} we have
\begin{equation*}
\begin{split}
| \langle \mathbb{E}_{\omega}[b,P_{k, \omega}]_1(f_1, f_2),f_3 \rangle |
& \le \E_{\omega} | \langle [b,P_{k,\omega}]_1(f_1, f_2),f_3 \rangle | \\
&\lesssim (1+\max k_i) \|f_1\|_{L^p(\R^{n+m})} \|f_2\|_{L^q(\R^{n+m})}
\|f_3\|_{L^{r'}(\R^{n+m})}. 
\end{split}
\end{equation*}

Our remaining task is to prove a  weak type estimate, which combined with the Banach range boundedness
implies Theorem \ref{thm:com1PP} via interpolation. The excellent general idea of using weak type estimates and interpolation in this spirit appears e.g. in the work of Muscalu--Pipher--Tao--Thiele \cite{MPTT}. That paper dealt with special singular integrals, namely bi-parameter bilinear multipliers. Such operators are paraproduct free. They were not considering
commutator estimates either.

Let $p,q \in (1, \infty)$ and $r \in (1/2,1)$ satisfy $1/p+1/q=1/r$.
We will show that given $f_1 \in L^p(\R^{n+m})$, $f_2 \in L^q(\R^{n+m})$ and a set $E \subset \R^{n+m}$ with $0 < |E| < \infty$,
there exists a subset $E' \subset E$ such that $|E'| \ge |E|/2$ and such that for all functions $f_3$
satisfying $|f_3| \le 1_{E'}$ there holds
\begin{equation}\label{eq:ResAvePP}
\begin{split}
| \langle \mathbb{E}_{\omega}&[b,P_{k,\omega}]_1(f_1, f_2),f_3 \rangle | \\
&\lesssim (1+\max k_i) \|f_1\|_{L^p(\R^{n+m})} \|f_2\|_{L^q(\R^{n+m})}|E|^{1/r'}.
\end{split}
\end{equation}

To prove \eqref{eq:ResAvePP} we consider the different types of partial paraproducts separately.
Here we assume that every $P_{k,\omega}$ is of the form \eqref{eq:ParticularPartialP}.
All the other types are handled with similar arguments and we will comment on this in the end of this proof.
The commutators are split using the identities from Section \ref{sec:marprod}. 
Define $P_{k,\omega}^b(f_1,f_2)$ to be equal to
$$
\sum_{\substack{K \in \calD^n_{\omega_1} \\ V \in \calD^m_{\omega_2}}} 
\sum_{\substack{I_1, I_2, I_3 \in \calD^n_{\omega_1} \\ I_i^{(k_i)} = K}} 
(\langle b \rangle_{I_3 \times V} - \langle b \rangle_{I_1 \times V}) a_{K, V, (I_i)}^{\omega} \bla  f_1, h_{I_1}^0 \otimes h_V  \bra 
 \Big\langle f_2, h_{I_2} \otimes \frac{1_V}{|V|} \Big\rangle 
h_{I_3}\otimes \frac{1_V}{|V|}.
$$
For an arbitrary $f_3$ we use \eqref{eq:1EX} to write that
\begin{equation}\label{eq:SplitComPartialP}
\begin{split}
&\langle [b,P_{k,\omega}]_1(f_1, f_2),f_3 \rangle
=  \sum_{i=1}^2\langle P_{k,\omega}(f_1, f_2),a_{i,\omega_1}^1(b,f_3) \rangle
-\sum_{i=1}^2\langle P_{k,\omega}(a_{i,\omega_2}^2(b,f_1), f_2),f_3 \rangle \\
&+  \sum_{\substack{K \in \calD^n_{\omega_1} \\ V \in \calD^m_{\omega_2} }} \sum_{\substack{I_1, I_2, I_3 \in \calD^n_{\omega_1} \\ I_i^{(k_i)} = K}}
a_{K, V, (I_i)}^{\omega} \bla  f_1, h_{I_1}^0 \otimes h_V  \bra 
 \Big\langle f_2, h_{I_2} \otimes \frac{1_V}{|V|} \Big\rangle 
\bla (\langle b \rangle_{I_3,1} - \langle b \rangle_{I_3 \times V}) \langle f_3, h_{I_3} \rangle_1 \bra_{V} \\
&-\sum_{\substack{K \in \calD^n_{\omega_1} \\ V \in \calD^m_{\omega_2}}} \sum_{\substack{I_1, I_2, I_3 \in \calD^n_{\omega_1} \\ I_i^{(k_i)} = K}} 
\bigg[a_{K, V, (I_i)}^{\omega} \bla (\langle b \rangle_{V,2} - \langle b \rangle_{I_1 \times V}) \langle f_1, h_{V} \rangle_2 , h_{I_1}^0\bra \\
& \hspace{4cm} \times  \Big\langle f_2, h_{I_2} \otimes \frac{1_V}{|V|} \Big\rangle 
\Big \langle f_3, h_{I_3}\otimes \frac{1_V}{|V|} \Big \rangle \bigg]  \\
&+\langle P^b_{k,\omega}(f_1,f_2),f_3 \rangle.
\end{split}
\end{equation}
Here we denoted $a_{i,\omega_1}^1(b,\cdot):=a_{i,\calD^n_{\omega_1}}^1(b,\cdot)$ and similarly with $a_{i,\omega_2}^2(b,\cdot)$.

Consider the term from the third line. 
We show that there exists a set $E' \subset E$ with $|E'| \ge \frac{99}{100}|E|$ so that for all $f_3$ such that $|f_3| \le 1_{E'}$ there holds that
\begin{equation}\label{eq:OneFromPartialPCom}
\begin{split}
\E_\omega \sum_{\substack{K \in \calD^n_{\omega_1} \\ V \in \calD^m_{\omega_2}}} \sum_{\substack{I_1, I_2, I_3 \in \calD^n_{\omega_1} \\ I_i^{(k_i)} = K}} 
&\Big| a_{K, V, (I_i)}^{\omega} \bla (\langle b \rangle_{V,2} - \langle b \rangle_{I_1 \times V}) 
\langle f_1, h_{V} \rangle_2 , h_{I_1}^0\bra \\
& \times  \Big\langle f_2, h_{I_2} \otimes \frac{1_V}{|V|} \Big\rangle 
\Big \langle f_3, h_{I_3}\otimes \frac{1_V}{|V|} \Big \rangle \Big| 
 \lesssim \| f_1 \|_{L^p} \| f_2 \|_{L^q} |E|^{1/r'}.
\end{split}
\end{equation}
The corresponding estimate for all the other terms can be proved with analogous arguments. We will briefly indicate the required modifications in the end of this proof.
Together these prove \eqref{eq:ResAvePP} for partial paraproducts that are of the form \eqref{eq:ParticularPartialP}.

Now, we turn to prove \eqref{eq:OneFromPartialPCom}. First, let  $\varphi^2_{\omega_2,b} f_1 := \varphi^2_{\calD^m_{\omega_2},b} f_1$ be the function from Section \ref{sec:AdapMaxFunc}.
We have
\begin{equation}\label{eq:phi_bDomination}
\big | \bla (\langle b \rangle_{V,2} - \langle b \rangle_{I_1 \times V}) 
\langle f_1, h_{V} \rangle_2 , h_{I_1}^0\bra \big |
\le \langle \varphi^2_{\omega_2,b} f_1, h^0_{I_1} \otimes h_V \rangle.
\end{equation}
Given $\omega_2$ define the square function $\wt S_{\omega_2}$ acting on functions $g \colon \R^m \to \C$ by
$$
\wt{S}_{\omega_2}g
= \Big(\sum_{V \in \calD^m_0} |\langle g, h_{V+\omega_2} \rangle|^2 \frac{1_V}{|V|} \Big)^{1/2},
$$ 
and then for $f \colon \R^{n+m} \to \C$ define
$$
\wt{S}^2_{\omega_2}f= \Big(\sum_{V \in \calD^m_0} |\langle f, h_{V+\omega_2} \rangle_2|^2 \otimes  \frac{1_V}{|V|} \Big)^{1/2}.
$$
Recall that $\E_{\omega}=\E_{\omega_1} \E_{\omega_2}$. Let  $\Phi_1$ and $\Phi_2^l$, $0 \le l \in \Z$, be the auxiliary operators
$$
\Phi_1(f)=\E_{\omega_2}M^1 \wt{S}^2_{\omega_2}( \varphi^2_{\omega_2,b} f), \quad 
\Phi_2^l(f)=\Big( \sum_{K \in \calD^n_0} \E_{\omega_1}(M^1\Delta^1_{K+\omega_1,l} \varphi^1_{\omega_1} f)^2 \Big)^{1/2},
$$
where $\varphi^1_{\omega_1}f:= \varphi^1_{\calD^n_{\omega_1}}f$ was introduced in Lemma \ref{lem:standardEst2}.

\begin{lem}\label{lem:AuxBdd}
We have for all $l \in \Z$, $l \ge 0$, that
$$
\|\Phi_1(f)\|_{L^s(\R^{n+m})}+\|\Phi_2^l(f)\|_{L^s(\R^{n+m})} \lesssim \|f\|_{L^s(\R^{n+m})}, \qquad s \in (1,\infty),
$$
where the bound is independent of $l$.
\end{lem}
\begin{proof}
We use weights and extrapolation (it is well-known that standard extrapolation results also work with bi-parameter weights). 
This is useful with $\Phi_2^l$ in order to reduce to $L^2$ estimates where we can take the expectation out. With $\Phi_1$ we could do without weights
by estimating directly in $L^s$. Take $w \in A_2(\R^n \times \R^m)$. 

Notice that for all $V \in \calD^m_0$ we have
$$
|\langle f, h_{V+\omega_2} \rangle_2| \otimes  \frac{1_V}{|V|^{1/2}} \lesssim M^2 ( \langle f, h_{V+\omega_2} \rangle_2 \otimes h_{V+\omega_2}).
$$
This implies
$$
\| \wt{S}^2_{\omega_2}f \|_{L^2(w)} \lesssim \Big\| \Big( \sum_{V \in \calD^m_{\omega_2}} [M^2 ( \langle f, h_{V} \rangle_2 \otimes h_{V})]^2 \Big)^{1/2} \Big\|_{L^2(w)}
\le C([w]_{A_2}) \|f\|_{L^2(w)},
$$
where we used weighted maximal function and weighted square function estimates. Lemma \ref{lem:bmaxbounds} says that
$\|\varphi_{\omega_2, b}^2(f)\|_{L^2(w)} \le C([w]_{A_2}) \|f\|_{L^2(w)}$, and obviously $M^1$ satisfies the same bound. The $L^2(w)$
result for $\Phi_1$ follows, and we can extrapolate.

Next, we have
\begin{align*}
\|\Phi_2^l(f)\|_{L^2(w)}^2 = \E_{\omega_1} \Big\| \Big( \sum_{K \in \calD^n_{\omega_1}} (M^1\Delta^1_{K,l} \varphi^1_{\omega_1} f)^2 \Big)^{1/2} \Big\|_{L^2(w)}^2
\le C([w]_{A_2}) \|f\|_{L^2(w)}^2
\end{align*}
using weighted maximal function and weighted square function estimates and Lemma \ref{lem:standardEst2}. 
We can extrapolate to finish.
\end{proof}

We have $\| \Phi_1(f_1)\Phi^{k_2}_2(f_2) \|_{L^r} \lesssim \| f_1 \|_{L^p} \| f_2 \|_{L^q}$. Recall Equation \eqref{eq:phi_bDomination}.
 Therefore, to prove \eqref{eq:OneFromPartialPCom} we may assume that
$\| \Phi_1(f_1)\Phi_2^{k_2}(f_2) \|_{L^r}=1$ and then
show that there exists a set $E' \subset E$ with $|E'| \ge \frac{99}{100}|E|$ so that for all $f_3$ such that $|f_3| \le 1_{E'}$ there holds that
\begin{equation}\label{eq:JustAsPartialP}
\begin{split}
\sum_{\substack{K \in \calD^n_0 \\ V \in \calD^m_0}} \Lambda_{K,V}(f_1,f_2,f_3)
\lesssim |E|^{1/r'},
\end{split}
\end{equation} 
where $\Lambda_{K,V}$ is defined to act on three functions by
\begin{equation*}
\begin{split}
\Lambda_{K,V}(f_1,f_2,f_3)=
\E_\omega \sum_{\substack{I_1, I_2, I_3 \in \calD^n_{\omega_1} \\ I_i^{(k_i)} = K+\omega_1}} 
\Big| & a_{K+\omega_1, V+\omega_2, (I_i)}^{\omega} 
\langle \varphi^2_{\omega_2,b} f_1, h^0_{I_1} \otimes h_{V+\omega_2} \rangle \\ 
&\times \Big\langle f_2, h_{I_2} \otimes \frac{1_{V+\omega_2}}{|V|} \Big\rangle 
\Big \langle f_3, h_{I_3}\otimes \frac{1_{V+\omega_2}}{|V|} \Big \rangle \Big|.
\end{split}
\end{equation*}

Define the sets
$$
\Omega_u =  \{ \Phi_1(f_1)\Phi_2^{k_2}(f_2) > C2^{-u} | E | ^{-1/r}  \}, \quad u \ge 0.
$$
For a small enough $c=c(n,m) \in (0,1)$ define the enlargement by 
$
\wt{\Omega}_u= \{M1_{\Omega_u} > c \}.
$
The set $E'$ is defined by $E'=E \setminus \wt{\Omega}_0$. 
By choosing the constant $C$ in the definition of the sets $\Omega_u$ to be large enough, we have $|E'| \ge \frac{99}{100} |E |$.
Then, let $\widehat{\calR}_u$ be the collection of rectangles
$$
\widehat{\calR}_u= \Big\{ R \in \calD_0 \colon
|R \cap \Omega_u| \ge \frac{1}{100} |R| \Big \},
$$
and write $\calR_u = \widehat{\calR}_u \setminus \widehat{\calR}_{u-1}$ when $u \ge 1$.

Now, we fix an arbitrary function $f_3$  such that $|f_3| \le 1_{E'}$ and consider \eqref{eq:JustAsPartialP} with $f_3$.
Let $K \times V \in \calD_0$. 
First, we show that if $\Lambda_{K,V} (f_1,f_2,f_3) \not=0$ then $K \times V \in \widehat \calR_u$ for some $u$.
We have for all $\omega=(\omega_1,\omega_2)$ and almost every $(x_1,x_2) \in K \times V$ that
\begin{equation*}
\begin{split}
\sum_{\substack{I_1 \in \calD^n_{\omega_1} \\ I_1^{(k_1)} = K+\omega_1}}
\frac{|I_1|^{1/2}}{|K||V|^{1/2}}| \langle \varphi^2_{\omega_2,b} f_1, h^0_{I_1} \otimes h_{V+\omega_2} \rangle|
& \le \sum_{\substack{I_1 \in \calD^n_{\omega_1} \\ I_1^{(k_1)} = K+\omega_1}} \frac{|I_1|^{1/2}}{|K|} \langle \wt S^2_{\omega_2} (\varphi^2_{\omega_2,b} f_1), h_{I_1}^0 \rangle_1(x_2) \\
&\lesssim M^1 \wt S^2_{\omega_2} (\varphi^2_{\omega_2,b} f_1)(x_1,x_2)
\end{split}
\end{equation*}
and
\begin{equation*}
\begin{split}
\sum_{\substack{I_2 \in \calD^n_{\omega_1} \\ I_2^{(k_2)} = K+\omega_1}} 
&\frac{|I_2|^{1/2}}{|K|} \Big | \Big \langle f_2, h_{I_2} \otimes \frac{1_{V+\omega_2}}{|V|} \Big \rangle \Big |
\lesssim 
 \sum_{\substack{I_2 \in \calD^n_{\omega_1} \\ I_2^{(k_2)} = K+\omega_1}} \frac{|I_2|^{1/2}}{|K|}
 M\langle f_2, h_{I_2} \rangle_1(x_2) \\
&= \sum_{\substack{I_2 \in \calD^n_{\omega_1} \\ I_2^{(k_2)} = K+\omega_1}} \frac{|I_2|^{1/2}}{|K|}
 \langle  \varphi^1_{\omega_1} f_2, h_{I_2} \rangle_1 (x_2)
\lesssim M^1(\Delta^1_{K+\omega_1,k_2} \varphi^1_{\omega_1} f_2) (x_1,x_2).
\end{split}
\end{equation*}
Therefore, if $\Lambda_{K,V} (f_1,f_2,f_3) \not=0$, then for
almost every $x \in K \times V$ there holds that
\begin{equation*}
\begin{split}
0 < \E_\omega \sum_{\substack{I_1,I_2 \in \calD^n_{\omega_1} \\ I_i^{(k_i)} = K+\omega_1}}
 \frac{|I_1|^{1/2}|I_2|^{1/2}}{|K|^2|V|^{1/2}} 
\Big| \langle \varphi^2_{\omega_2,b} f_1, h^0_{I_1} & \otimes h_{V+\omega_2} \rangle
 \Big \langle f_2, h_{I_2} \otimes \frac{1_{V+\omega_2}}{|V|} \Big \rangle \Big | \\
& \lesssim \Phi_1(f_1)(x) \Phi_2^{k_2}(f_2)(x).
\end{split}
\end{equation*}
The inequality ``$<$'' holds since the integrand is positive for $\omega$ in a set of positive measure.
From this it follows that $K \times V \subset \Omega_u$ if $u$ is large enough, and so in particular $K \times V \in \widehat \calR_u$.

If $K \times V \in \widehat{\calR}_u$, then for all $\omega$ there holds that $(K \times V)+\omega \subset (3K) \times (3V) \subset \wt \Omega_u$.
The constant $c=c(n,m)$ in the definition of $\wt \Omega_u$ is chosen so that this inclusion holds.
If $K \times V \in \widehat \calR _0$, then $(K \times V)+\omega \subset \wt \Omega_0 \subset (E')^c$ for all $\omega$, which combined with the fact that
$|f_3| \le 1_{E'}$ implies that $\Lambda_{K,V}(f_1,f_2,f_3)=0$.

With the above observations we have that
$$
\sum_{\substack{K \in \calD^n_0 \\ V \in \calD^m_0}} \Lambda_{K,V}(f_1,f_2,f_3)
= \sum_{u=1}^\infty \sum_{K \times V \in \calR_u} 
\Lambda_{K,V}(f_1,f_2,f_3).
$$
We fix $u$ and estimate the corresponding term.

First, notice that if $K \times V \in \calR_u$, then 
$\Lambda_{K,V}(f_1,f_2,f_3)=\Lambda_{K,V}( f_1,f_2,1_{\wt \Omega_u}f_3)$.
For $K \in \calD^n_0$ define
$
\calC_{K,u}
=\{V \in \calD^m_0 \colon K \times V \in \calR_u\}
$,
which allows us to write $\sum_{K \times V \in \calR_u}=\sum_{K \in \calD^n_0} \sum_{V \in \calC_{K,u}}$.
If $K \in \calD^n_0$, then every $V \in \calC_{K,u}$ satisfies $|(K \times V) \cap \Omega_{u-1}^c | \ge \frac{99}{100} |K \times V|$.
Therefore, Remark \ref{lem:H1-BMO-Modified-1par} gives that
\begin{equation}\label{eq:ApplyH1BMOLemma}
\begin{split}
& \sum_{K \in \calD^n_0}  \sum_{V \in \calC_{K,u}}\Lambda_{K,V}( f_1,f_2,1_{\wt \Omega_u}f_3) 
 \lesssim \sum_{K \in \calD^n_0} \E_\omega \sum_{\substack{I_1, I_2, I_3 \in \calD^n_{\omega_1} \\ I_i^{(k_i)} = K+\omega_1}} 
 \bigg[\frac{ |I_1|^{1/2}|I_2|^{1/2}|I_3|^{1/2}}{|K|^2} \\
&\hspace{1cm}\times \iint \displaylimits_{\wt\Omega_u \setminus \Omega_{u-1}}
\frac{1_{K}}{|K|} \otimes \wt S_{\omega_2} \langle \varphi^2_{\omega_2,b} f_1,  h^0_{I_1}  \rangle_1 M \langle f_2, h_{I_2} \rangle_1 M \langle  1_{\wt \Omega_u}f_3, h_{I_3} \rangle_1 \bigg], 
\end{split}
\end{equation}
where  we used the estimate 
$$
\Big | \Big \langle f_2, h_{I_2} \otimes \frac{1_{V+\omega_2}}{|V|} \Big \rangle \Big | 
\lesssim M \langle f_2, h_{I_2} \rangle_1(x_2), \quad x_2 \in V,
$$
and the same estimate with $1_{\wt \Omega_u}f_3$. 
We were able to insert the restriction $\wt \Omega_u$ to the integration area since $K \times V \subset \wt \Omega_u$ for every $K \times V \in \calR_u$.

Notice that 
$\wt S_{\omega_2} \langle \varphi^2_{\omega_2,b} f_1,  h^0_{I_1}  \rangle_1(x_2) 
\le \langle \wt S_{\omega_2}^2  (\varphi^2_{\omega_2,b} f_1),  h^0_{I_1}  \rangle_1(x_2)$ 
for all $x_2$, and 
recall that $M \langle f_2, h_{I_2} \rangle_1=\langle  \varphi^1_{\omega_1} f_2, h_{I_2} \rangle_1$. 
Thus, the inner sum over the cubes $I_i$ in the right hand side of \eqref{eq:ApplyH1BMOLemma} 
 is dominated by
\begin{equation*}
\begin{split}
\iint \displaylimits_{\wt\Omega_u \setminus \Omega_{u-1}}
M^1 \wt S^2_{\omega_2 }(\varphi^2_{\omega_2,b} f_1)
M^1(\Delta^1_{K+\omega_1,k_2} \varphi^1_{\omega_1} f_2)
M^1(\Delta^1_{K+\omega_1,k_3} \varphi^1_{\omega_1} (1_{\wt \Omega_u}f_3)).
\end{split}
\end{equation*}
Taking expectation $\E_{\omega}=\E_{\omega_1} \E_{\omega_2}$, using H\"older's inequality with respect to $\omega_1$ 
and summing over $K \in \calD^n_0$ shows that 
\begin{equation*}
\begin{split}
\sum_{K \in \calD^n_0}  \sum_{V \in \calC_{K,u}}\Lambda_{K,V}( f_1,f_2,1_{\wt \Omega_u}f_3) 
\lesssim \iint \displaylimits_{\wt\Omega_u \setminus \Omega_{u-1}} \Phi_1(f_1) \Phi_2^{k_2}(f_2)
\Phi_2^{k_3}(1_{\wt \Omega_u} f_3).
\end{split}
\end{equation*}

By definition we have $\Phi_1(f_1)(x) \Phi_2^{k_2}(f_2)(x) \lesssim 2^{-u}|E|^{-1/r}$ for all $x \in \Omega_{u-1}^c$. Also, just by using the 
$L^2$-boundedness of $\Phi_2^{k_3}$ and the fact that $\| f_3\|_{L^\infty} \le 1$ there holds that
$$
\iint \displaylimits_{\wt\Omega_u \setminus \Omega_{u-1}}
\Phi_2^{k_3}(1_{\wt \Omega_u} f_3)
\lesssim |\wt\Omega_u| \lesssim |\Omega_u| \lesssim 2^{ur} |E|.
$$
These combined give that
$$
\sum_{K \in \calD^n_0}  \sum_{V \in \calC_{K,u}}
\Lambda_{K,V}(f_1,f_2,f_3)  \lesssim 2^{-u(1-r)} |E|^{1/r'},
$$
which can be summed over $u$ since $r < 1$. This finishes the proof of \eqref{eq:JustAsPartialP}.

Let us now briefly comment on how to handle the other terms from \eqref{eq:SplitComPartialP}.
The main difference in the beginning is how to define the sets $\Omega_u$. After that, one proceeds with the corresponding
steps as above. For instance, when considering  $\langle P_{k,\omega}(a_{i,\omega_2}^2(b,f_1), f_2),f_3 \rangle$, one sets
$$
\Omega_u =  \{ \E_{\omega_2}M^1 \wt{S}^2_{\omega_2}( a^2_{i,\omega_2} f)\Phi_2^{k_2}(f_2) > C2^{-u} | E | ^{-1/r}  \}, \quad u \ge 0,
$$
where $\Phi_2^{k_2}$ is as above. Just like $\Phi_1$ in Lemma \ref{lem:AuxBdd} the 
operator $ f \mapsto \E_{\omega_2}M^1 \wt{S}^2_{\omega_2}( a^2_{i,\omega_2} f)$ is bounded.

The function $f_3$ has a special role in the proof, and for example localisation properties in the spirit of
$\Lambda_{K,V}(f_1,f_2,f_3)=\Lambda_{K,V}( f_1,f_2,1_{\wt \Omega_u}f_3)$ for $K \times V \in \calR_u$ are important.
One has to be careful with this regarding the other terms in \eqref{eq:SplitComPartialP}, and we deal with a different type of term
with the full paraproducts to make this even clearer.
For example, with the term $\langle P_{k,\omega}(f_1, f_2),a_{i,\omega_1}^1(b,f_3) \rangle$ one uses the fact that the operators
$a_{i,\omega_1}^1(b,f_3)$ have the  localisation property \eqref{eq:Localisation}.
The corresponding localisation is also important with the term from the second line of \eqref{eq:SplitComPartialP}.
Therefore, when dealing with this term, one must not use the domination as in \eqref{eq:phi_bDomination} right in the beginning, as we did above, 
but after the localisation is used.

Finally, we discuss the term $\langle P^b_{k,\omega}(f_1,f_2),f_3 \rangle$, which is in fact the easiest one (although in the linear Bloom
setting, see e.g. \cite{LMV3}, this is probably the hardest term).
Simply begin by using the estimate
\begin{equation*}
|\langle b \rangle_{I_3 \times V} - \langle b \rangle_{I_1 \times V}|
\le |\langle b \rangle_{I_3 \times V} - \langle b \rangle_{K \times V}|
+|\langle b \rangle_{K \times V} - \langle b \rangle_{I_1 \times V}|
\lesssim \| b \|_{\bmo(\R^{n+m})} \max_i k_i,
\end{equation*}
and then proceed in the usual way. This term produces the complexity dependency.
We have now commented on all the terms from \eqref{eq:SplitComPartialP}, which concludes the
proof for partial paraproducts that are of the form  \eqref{eq:ParticularPartialP}.

There is no essential difference with other types of partial paraproducts. 
First, one uses the relevant identities from Section \ref{sec:marprod}, which leads to a splitting analogous to 
\eqref{eq:SplitComPartialP}. 
To each type of partial paraproduct there is the natural combination of maximal functions and square functions, which are used 
to build the auxiliary operators that correspond to $\Phi_1$ and $\Phi^l_2$ above.
With the auxiliary operators the argument goes as before. 

We just remark the following. The operators $A_{i,\omega}(b, \cdot)$  from Section \ref{sec:marprod} appear
related to some types of partial paraproducts. All these satisfy the localisation property corresponding to  \eqref{eq:Localisation},
and auxiliary operators involving $A_{i,\omega}(b, \cdot)$ are bounded with the same proof as in Lemma \ref{lem:AuxBdd}.
Moreover, some terms involve pairings like $\langle  (b-\langle b \rangle_{I \times V}) f, h^0_I \otimes \frac{1_V}{|V|} \rangle$. These
clearly have the property analogous to \eqref{eq:Localisation}, and at some point of the argument one uses the estimate
$$
\Big | \Big \langle  (b-\langle b \rangle_{I \times V}) f, h^0_I \otimes \frac{1_V}{|V|} \Big  \rangle \Big |
\le \Big \langle M_b f , h^0_I \otimes \frac{1_V}{|V|} \Big  \rangle.
$$
Here $M_b$ is the maximal function from Section \ref{sec:AdapMaxFunc}. The operator $M_b$ is bounded as stated in Lemma \ref{lem:bmaxbounds}.
This ends the proof of Theorem \ref{thm:com1PP}.

\end{proof}

\section{Commutators of full paraproducts}
We begin by defining full paraproducts on random dyadic grids. Let $\omega = (\omega_1, \omega_2) \in (\{0,1\}^n)^{\Z} \times (\{0,1\}^m)^{\Z}$.
A full paraproduct $\Pi_{\omega}$ of a particular form is
$$
\langle \Pi_{\omega}(f_1, f_2), f_3 \rangle = \sum_{K \times V \in \calD_{\omega}} a_{K,V}^{\omega} \langle f_1 \rangle_{K \times V} \langle f_2 \rangle_{K \times V} 
\langle f_3, h_K \otimes h_V\rangle,
$$
where
$$
\|\{a_{K,V}^{\omega}\}_{K \times V \in \calD_{\omega}}\|_{\BMO_{\textup{prod}}(\calD_{\omega})} \le 1.
$$
There are again nine different forms of full paraproducts: if we view
$\langle f_1 \rangle_{K \times V} = \big\langle f_1, \frac{1_K}{|K|} \otimes \frac{1_V}{|V|} \big \rangle$ and similarly with $f_2$, we can put the the cancellative
Haar function $h_K$ (currently paired with $f_3$) to any of the other two slots that currently have $1_K/|K|$, and similarly with $h_V$.

\begin{thm}\label{thm:com1FP}
Let $\|b\|_{\bmo(\R^{n+m})}  = 1$, and let $1 < p, q \le \infty$ and $1/2 < r < \infty$ satisfy $1/p+1/q = 1/r$.
Suppose that $(\Pi_{\omega})_{\omega}$ is a collection of full paraproducts of the same type.
Then we have
$$
 \|\mathbb{E}_{\omega}[b,\Pi_{\omega}]_1(f_1, f_2)\|_{L^r(\R^{n+m})} 
 \lesssim \|f_1\|_{L^p(\R^{n+m})} \|f_2\|_{L^q(\R^{n+m})}.
$$
\end{thm}
\begin{proof}
As in the partial paraproduct case we are done after showing the following.
Let $p,q \in (1, \infty)$ and $r \in (1/2,1)$ satisfy $1/p+1/q=1/r$.
We will show that given $f_1 \in L^p(\R^{n+m})$, $f_2 \in L^q(\R^{n+m})$ and a set $E \subset \R^{n+m}$ with $0 < |E| < \infty$,
there exists a subset $E' \subset E$ such that $|E'| \ge |E|/2$ and such that for all functions $f_3$
satisfying $|f_3| \le 1_{E'}$ there holds
\begin{equation}\label{eq:ResAveFP}
| \langle \mathbb{E}_{\omega}[b,\Pi_{\omega}]_1(f_1, f_2),f_3 \rangle |
\lesssim  \|f_1\|_{L^p(\R^{n+m})} \|f_2\|_{L^q(\R^{n+m})}|E|^{1/r'}.
\end{equation}

As with partial paraproducts in Section \ref{sec:ComPartialP},  the different forms of full paraproducts are handled separately, 
but with analogous arguments. Let us consider here the case where every $\Pi_\omega$ is of the form
\begin{equation}\label{eq:OneFullPara}
\langle \Pi_{\omega}(f_1, f_2), f_3 \rangle = \sum_{K \times V \in \calD_{\omega}} 
a_{K,V}^{\omega} \langle f_1 \rangle_{K \times V} 
\Big \langle f_2, \frac{1_K}{|K|} \otimes h_V \Big \rangle  
\Big \langle f_3, h_K \otimes \frac{1_V}{|V|} \Big \rangle.  
\end{equation}
The commutators $[b,\Pi_\omega]$ are again split with the identities from Section \ref{sec:marprod}, this time
using \eqref{eq:1EX} and \eqref{eq:noEX}.
The resulting terms are handled separately with similar arguments. 
Here we consider the term
\begin{equation}\label{eq:OneFromComFullP}
\E_\omega \sum_{K \times V \in \calD_{\omega}} 
a_{K,V}^{\omega} \langle f_1 \rangle_{K \times V} 
\Big \langle f_2, \frac{1_K}{|K|} \otimes h_V \Big \rangle  
\Big \langle a^1_{i,\omega_1}(b, f_3), h_K \otimes \frac{1_V}{|V|} \Big \rangle
\end{equation}
for some $i \in \{1,2\}$. 
We will not discuss the other terms and the other forms of full paraproducts, but
just refer to the corresponding discussion in the end of the proof of Theorem \ref{thm:com1PP}.

This time, let $\Phi_1$ and $\Phi_2$ be the auxiliary operators
$$
\Phi_1(f)
=\E_{\omega_2}\Big( \sum_{V \in \calD^m_0} (M \langle f, h_{V+\omega_2} \rangle_2)^2 \otimes \frac{1_V}{|V|} \Big)^{1/2}
$$ 
and 
$$
\Phi_2(f)
=\E_{\omega_1}\Big( \sum_{K \in \calD^n_0} \frac{1_K}{|K|} \otimes (M \langle a^1_{i,\omega_1}(b,f), h_{K+\omega_1} \rangle_1)^2 \Big)^{1/2}.
$$ 
Similarly as in Lemma \ref{lem:AuxBdd} (notice that $\Phi_1$ and $\Phi_2$ have different definitions than in Lemma \ref{lem:AuxBdd}) 
these are bounded in $L^s$ for every $s \in (1, \infty)$, which uses the fact that 
$a^1_{i,\omega_1}(b,\cdot)$ is bounded by Lemma \ref{lem:basicAa}.

Let now $f_1 \in L^p(\R^{n+m})$ and $f_2 \in L^q(\R^{n+m})$ be such that  $\| Mf_1 \Phi_1(f_2) \|_{L^r}=1$ and let $E \subset \R^{n+m}$ with $0 < |E| < \infty$.
We show that there exists a set $E' \subset E$ with $|E'| \ge \frac{99}{100}|E|$ so that 
\begin{equation}\label{eq:EstOneFromComFullP}
 \sum_{K \times V \in \calD_{0}} \Lambda_{K,V}(f_1,f_2,f_3)
 \lesssim |E|^{1/r'}
\end{equation}
holds for all $f_3$ such that $|f_3| \le 1_{E'}$,
where $\Lambda_{K,V}$ acts on a triple of functions by
\begin{equation*}
\begin{split}
\Lambda_{K,V}(f_1,f_2,f_3)
=\E_\omega
\Big |a_{K+\omega_1,V+\omega_2}^{\omega} \langle f_1 \rangle_{(K \times V)+\omega} 
&\Big \langle f_2, \frac{1_{K+\omega_1}}{|K|} \otimes h_{V+\omega_2} \Big \rangle  \\
&\times \Big \langle a^1_{i,\omega_1}(b, f_3), h_{K+\omega_1} \otimes \frac{1_{V+\omega_2}}{|V|} \Big \rangle \Big |.
\end{split}
\end{equation*}
Since $\| Mg_1 \Phi_1(g_2)\|_{L^r} \lesssim \| g_1 \|_{L^p} \| g_2 \|_{L^q}$ for all $g_1$ and $g_2$, this gives the estimate that we want for
the term \eqref{eq:OneFromComFullP}. Together with the corresponding estimates for all the other parts of 
$\E_\omega[b,\Pi_\omega]$, this proves \eqref{eq:ResAveFP} for full paraproducts that are of the form \eqref{eq:OneFullPara}.

We turn to prove \eqref{eq:EstOneFromComFullP}. For $u \ge 0$
let $\Omega_u=\{Mf_1\Phi_1(f_2) > C2^{-u} |E|^{-1/r}\}$ and
$\wt \Omega_u =\{M1_{\Omega_u} > c\}$, where $c=c(n,m) \in (0,1)$ is a small constant. 
Set $E' = E \setminus \wt \Omega_0$.  By choosing the constant $C$ to be large enough,
we have that $|E'| \ge \frac{99}{100}|E|$.
Then, define $\widehat{\calR}_u= \{R \in \calD_0 \colon |R \cap \Omega_u| \ge \frac{1}{100}|R|\}$ for $u \ge 0 $ and
$\calR_u=\widehat{\calR}_u \setminus \widehat{\calR}_{u-1}$ for $u \ge 1$.

Suppose $K \times V \in \calD_0$ is such that $\Lambda_{K,V} (f_1,f_2,f_3) \not=0$.
Then, wee see that
$$
0 <  \E_\omega \frac{1}{|V|^{1/2}} 
\Big |\langle f_1 \rangle_{(K \times V)+\omega} 
\Big \langle f_2, \frac{1_{K+\omega_1}}{|K|} \otimes h_{V+\omega_2} \Big \rangle \Big |
\lesssim Mf_1(x)\Phi_1(f_2)(x)
$$
for all $x \in K \times V $. The first ``$<$'' holds since the integrand is positive for $\omega$ in a set of positive measure.
Thus, $K \times V \subset \Omega_u$ for large enough $u$, so $K \times V \in \widehat{\calR}_u$.

If $R \in \widehat{\calR}_u$, then $R +\omega \subset 3R \subset  \wt \Omega_u$ for all $\omega$, which is based the fact that $c=c(n,m)$ in the definition of $\wt \Omega_u$
is small enough. Notice that 
\begin{equation}\label{eq:Localisation}
\Big \langle a^1_{i,\omega_1}(b, f_3), h_{K+\omega_1} \otimes \frac{1_{V+\omega_2}}{|V|} \Big \rangle
=\Big \langle a^1_{i,\omega_1}(b, 1_{(K\times V)+\omega }f_3), h_{K+\omega_1} \otimes \frac{1_{V+\omega_2}}{|V|} \Big \rangle.
\end{equation}
Thus, if $K \times V \in \widehat{\calR}_0$, then $\Lambda_{K,V}(f_1,f_2,f_3)=0$ since $(K \times V)+\omega \subset (E')^c$ for all $\omega$, and $|f_3| \le 1_{E'}$.

Now, we have that
$$
\sum_{K \times V \in \calD_{0}} \Lambda_{K,V}(f_1,f_2,f_3)
= \sum_{u=1}^\infty \sum_{K \times V \in \calR_u} \Lambda_{K,V}(f_1,f_2,1_{ \wt \Omega_u}f_3),
$$
where it was legitimate to replace $f_3$ with $1_{ \wt \Omega_u}f_3$ because of \eqref{eq:Localisation}. We fix one $u$ and 
estimate the related term.

If $K \times V \in \calR_u$, then by definition $|(K \times V) \cap \Omega_{u-1}^c| \ge \frac{99}{100} |K \times V|$.
Therefore, using Lemma \ref{lem:H1-BMO-Modified} and then
the estimate
$$
\Big |\Big \langle f_2, \frac{1_{K+\omega_1}}{|K|} \otimes h_{V+\omega_2 } \Big\rangle \Big |
\lesssim M \langle f_2, h_{V+\omega_2} \rangle_2 (x_1), \quad x_1 \in K,
$$
and a corresponding estimate related to $f_3$, we have that
\begin{equation*}
\begin{split}
\sum_{K \times V \in \calR_u} \Lambda_{K,V}(f_1,f_2,1_{\wt \Omega_u} f_3) 
\lesssim  \iint   \displaylimits_{\wt \Omega_u \setminus \Omega_{u-1}}  Mf_1 \Phi_1(f_2) \Phi_2(1_{\wt \Omega_u} f_3).
\end{split}
\end{equation*}
The restriction to $\wt \Omega_u$ in the integration came from the fact that every $R \in \calR_u$ satisfies $R \subset \wt \Omega_u$.
Since $Mf_1(x) \Phi_1(f_2)(x) \lesssim 2^{-u}|E|^{-1/r}$ for $x \in \Omega_{u-1}^c$, the operator $\Phi_2$ is $L^2$ bounded and $\|f_3\|_{L^\infty} \le 1$,
there holds that 
$$
\iint   \displaylimits_{\wt \Omega_u \setminus \Omega_{u-1}}  Mf_1 \Phi_1(f_2) \Phi_2(1_{\wt \Omega_u} f_3)
\lesssim 2^{-u}|E|^{-1/r} |\wt \Omega_u|
\lesssim 2^{-u(1-r)}|E|^{1/r'}.
$$
This can be summed over $u$ since $r < 1$, which finishes the proof of \eqref{eq:EstOneFromComFullP}.
\end{proof}

\section{Synthesis: Proof of Theorem \ref{thm:main}} 
We can now use the representation theorem valid for bilinear biparameter singular integrals \cite[Theorem 5.1]{LMV1}, the boundedness of commutators of shifts
in the full range \cite{LMV2}, and the boundedness of commutators of partial and full paraproproducts in the full range -- proved in the current paper -- to give a short proof
of Theorem \ref{thm:main} 

\begin{proof}[Proof of Theorem \ref{thm:main}]
Using \cite[Theorem 5.1]{LMV1} we write the pointwise identity
$$
[b,T]_1(f_1, f_2) = C_T \mathop{\sum_{k = (k_1, k_2, k_3) \in \Z_+^3}}_{v = (v_1, v_2, v_3) \in \Z_+^3} \alpha_{k, v} 
\sum_{u} \mathbb{E}_{\omega} [b,U^{v}_{k, u, \mathcal{D}_{\omega}}]_1(f_1, f_2),
$$
where $C_T \lesssim 1$, $\alpha_{k, v} = 2^{- \alpha \max k_i/2} 2^{- \alpha \max v_j/2}$ ($\alpha > 0$ appears in the kernel estimates of $T$),
the summation over $u$ is finite, and
$U^{v}_{k, u, \mathcal{D}_{\omega}}$ is always either a dyadic shift of complexity $(k,v)$,
a partial paraproduct of complexity $k$ or $v$ (this requires $k= 0$ or $v=0$) or a full paraproduct (this requires $k=v=0$). 

We use $\|\sum_i g_i\|_{L^r}^r \le \sum_i \| g_i  \|_{L^r}^r$ if $r < 1$, and otherwise we use the normal triangle inequality. Using then that always
$$
 \|\mathbb{E}_{\omega}[b,U^{v}_{k, u, \mathcal{D}_{\omega}}]_1(f_1, f_2)\|_{L^r(\R^{n+m})} 
 \lesssim (1+\max(k_i, v_i)) \|f_1\|_{L^p(\R^{n+m})} \|f_2\|_{L^q(\R^{n+m})}
$$
we get the claim. For shifts this is proved in \cite{LMV2}, and for partial and full paraproducts this is proved above.
Notice that it was critical to consider averages of model operators in the range $r < 1$.
\end{proof}

\end{document}